\documentclass[12pt]{amsart}
\usepackage{amscd,amssymb,amsmath,latexsym,graphicx}

\newtheorem{thm}{Theorem}[section]
\newtheorem{lem}[thm]{Lemma}
\newtheorem{prop}[thm]{Proposition}
\newtheorem{cor}[thm]{Corollary}

\newtheorem{rmk}{Remark}[section]
\newtheorem{ex}[rmk]{Example}

\newcommand \SnG{\text{$\Sigma^n(G)$}}
\newcommand \OnG{\text{$\Omega^n(G)$}}
\newcommand \OnH{\text{$\Omega^n(H)$}}
\newcommand \OnGH{\text{$\Omega^n(G \times H)$}}
\newcommand \ep{\text{$e^{\prime}$}}
\newcommand \epp{\text{$e^{\prime\prime}$}}

\newcommand \Rm{\text{$\mathbb{R}^m$}}

\newcommand \Hom{\text{Hom$(G,\mathbb{R})$ }}
\newcommand \bdR{\text{$\partial_{\infty}\Rm$}}

\newcommand \hgs{\text{$H_{\gamma,s}$}}
\newcommand \cgs{\text{$C_{\gamma,s}$}}
\newcommand \xgs{\text{$X_{\gamma,s}$}}
\newcommand \xgsl{\text{$X_{\gamma,s-\lambda}$}}
\newcommand \ygs{\text{$Y_{\gamma,s}$}}
\newcommand \ygsl{\text{$Y_{\gamma,s-\lambda}$}}

\begin{document}
\pagestyle{plain}

\title{A Relationship Between Twisted Conjugacy Classes and the Geometric Invariants $\Omega^n$}
\author{Nic Koban and Peter Wong}
\date{\today}
\thanks{The second author was supported in part by the National Science Foundation DMS-0805968.}
\maketitle

\begin{abstract}
A group $G$ is said to have the property $R_\infty$ if every automorphism $\varphi \in {\rm Aut}(G)$ has an infinite number of $\varphi$-twisted conjugacy classes. Recent work of Gon\c calves and Kochloukova uses the $\Sigma^n$ (Bieri-Neumann-Strebel-Renz) invariants to show the $R_{\infty}$ property for a certain class of groups, including the generalized Thompson's groups $F_{n,0}$. In this paper, we make use of the $\Omega^n$ invariants, analogous to $\Sigma^n$, to show $R_{\infty}$ for certain finitely generated groups. In particular, we give an alternate and simpler proof of the $R_{\infty}$ property for $BS(1,n)$. Moreover, we give examples for which the $\Omega^n$ invariants can be used to determine the $R_{\infty}$ property while the $\Sigma^n$ invariants techniques cannot.
\end{abstract}

\noindent
{\bf 2000 Mathematics Subject Classification:} primary 20F65, 20E45; secondary 55M20, 57M07

\noindent
{\bf Keywords:} twisted conjugacy class, property $R_{\infty}$, $\Sigma$-invariants, $\Omega$-invariants

\section{Introduction}

An old problem in group theory asks whether an infinite group must have an infinite number of conjugacy classes of elements. In 1949 \cite{HNN}, Higman-Neumann-Neumann constructed an infinitely generated group with only a finite number of conjugacy classes. Ivanov (see \cite{I1,I2}) announced without proof an example of a finitely generated infinite group with a finite number of conjugacy classes. More recently, D. Osin \cite{osin} constructed a finitely generated infinite group in which any two non-trivial elements are conjugate. A more general question is: if $G$ has an infinite number of conjugacy classes and $\varphi$ is an automorphism of $G$, then does there exist an infinite number of $\varphi$-twisted conjugacy classes? Denote by $R(\varphi)$ the number of $\varphi$-twisted conjugacy classes or equivalently the number of orbits of the (left) action of $G$ on $G$ via $\sigma \cdot \alpha \mapsto \sigma \alpha \varphi(\sigma)^{-1}$. Then we can ask whether $R(\varphi)=\infty$. For instance, when $G=\mathbb Z$, the automorphism group is ${\rm Aut}(\mathbb Z)=\{1_{\mathbb Z}, -1_{\mathbb Z}\}$. It is easy to see that $R(1_{\mathbb Z})=\infty$ whereas $R(-1_{\mathbb Z})=2$.

The classical Burnside-Frobenius theory states that the number of equivalence classes of irreducible unitary representations of a finite group $G$ is equal to the number of conjugacy classes of elements of $G$. More generally, if $\varphi$ is an automorphism of $G$, the number of fixed points for the action of $\varphi$ on the unitary dual $\widehat G$ is equal to $R(\varphi)$ when $G$ is finite \cite{FH}. For infinite groups, this equality does not hold in general. Fel'shtyn and Troitsky \cite{FT} showed that this equality does hold for virtually polycyclic groups if one of the numbers in the equality is finite.

The number $R(\varphi)$ also plays an important role in classical Nielsen fixed point theory (for more details, see for example \cite{jiang}). Given a group endomorphism $\varphi:G \to G$, the $\varphi$-twisted conjugacy classes are in one-one correspondence with the fixed point (Reidemeister) classes (including the empty ones via the covering space approach \cite{jiang}) of a continuous map $f:M\to M$ on a compact connected polyhedron $M$ with $G=\pi_1(M)$ and $\varphi =f_{\sharp}$ the induced homomorphism on $\pi_1$. The Reidemeister number $R(f)$ is defined to be the cardinality of the set of the $\varphi=f_{\sharp}$ twisted conjugacy classes. It is well defined and is independent of the basepoint. The number $R(f)$ is always an upper bound for the Nielsen number $N(f)$ and for a large class of spaces, which include the so-called Jiang spaces, the Nielsen number $N(f)$ is either $0$ or is equal to $R(f)$. Since the Nielsen number is often equal to the minimum number of fixed points in the homotopy class of $f$, the computation of $N(f)$ is a central issue but is very difficult in general. Therefore, a Jiang-type result where $N(f)=0$ or $N(f)=R(f)$ is very desirable. In particular, if such a space is a manifold of dimension at least $3$, $R(f)=\infty$ implies that $f$ is deformable to be fixed point free. Furthermore, for a homeomorphism $f$ on such manifolds of dimension at least $5$, $R(f)=\infty$ implies that $f$ is isotopic to a fixed point free homeomorphism. Therefore, if $R(\varphi)=\infty$ for every automorphism $\varphi \in {\rm Aut}(\pi_1(M))$ then {\it every} homeomorphism can be made fixed point free. For instance in \cite{daci-peter3}, we have constructed for any integer $n\ge 5$ compact nilmanifolds and compact flat manifolds (or equivalently infra-abelian manifolds) of dimension $n$ such that every homeomorphism is isotopic to be fixed point free. It would be interesting to find examples of groups $G$ where $R(\varphi)=\infty$ {\it for every} $\varphi \in {\rm Aut}(G)$. In fact, Fel'shtyn and Hill conjectured \cite{FH} that $R(\varphi)=\infty$ for every injective endomorphism of a finitely generated group with exponential growth. This was later shown to be false (see \cite{daci-peter1}) in general.

We have seen in both group theory and in fixed point theory that the finiteness of $R(\varphi)$ is of great interest especially when $\varphi$ is an automorphism.
Following the same terminology in \cite{TW1}, we say that a group $G$ has the property $R_{\infty}$ if for every automorphism $\varphi \in {\rm Aut}(G)$, $R(\varphi)=\infty$. In recent years, there has been a growing interest in the study of the $R_{\infty}$ property. Finitely generated groups that possess this property include non-elemetary Gromov hyperbolic groups \cite{levitt-lustig}(see also \cite{F2}); non-elementary GBS (Generalized Baumslag Solitar) groups \cite{levitt} and groups that are quasi-isometric to such GBS groups \cite{TW2}; Baumslag-Solitar groups $BS(m,n)$ for $(m,n)\neq (1,1)$ \cite{fel-daci}; certain solvable generalization of $BS(1,n)$ and groups that are quasi-isometric to such groups \cite{TW1}; lamplighter groups $L_n$ where $gcd(n,6)\neq 1$ \cite{daci-peter2}, the Thompson's group $F$ \cite{BFG}, certain nilpotent groups \cite{daci-peter3}, certain polycyclic groups \cite{daci-peter1}, most of the 17 wallpaper groups \cite{daci-peter4}, a wide class of saturated weakly branch groups, which includes the Grigorchuk group of intermediate growth and the Gupta-Sidki group \cite{FLT}. Both combinatorial and geometric group theoretic methods have been employed in establishing the $R_{\infty}$ property for various families of groups.

More recently, the $\Sigma^n$ invariants of Bieri-Neumann-Strebel-Renz were used in \cite{GK} to show that the generalized Thompson's groups $F_{n,0}$, among other families of groups, have the $R_{\infty}$ property. The main idea in \cite{GK} is to show that if there is a finite set of discrete points on the character sphere that is ${\rm Aut}(G)$-invariant then $G$ has the $R_{\infty}$ property. Gon\c calves and Kochloukova \cite{GK} were able to identify certain groups $G$ with such property by computing the invariants $\Sigma^1(G)$. In general, the computation of $\Sigma^n(G)$ is difficult and very few such invariants have been computed. On the other hand, the $\Omega^n$ invariants, which seem easier to compute than $\Sigma^n$, can be computed for direct products (see \cite{K2}). The main objective of this paper is to employ the $\Omega^n$ invariants to study the $R_{\infty}$ property. Our results complement those obtained in \cite{GK}. We also give examples where our method may determine the $R_{\infty}$ property while the approach of \cite{GK} is not applicable.

The invariants $\Omega^n$ were defined in \cite{K1} and are analogs of the Bieri-Neumann-Strebel-Renz invariants $\Sigma^n$ defined in \cite{BNS} for $n=1$ and in \cite{BR} for $n \geq 2$.  We recall these definitions here for $n=1$; the full definitions are given in \S~\ref{sigma} and \S~\ref{omega}.

Let $G$ be a finitely generated group with generating set $\mathcal{X}$.  The set \Hom of homomorphisms from $G$ to the additive group of reals is a real vector space with dimension equal to the ${\mathbb Z}$-rank of the abelianization of $G$, so \Hom $\cong \Rm$ for some $m$. Choose an inner product in \Hom .  Denote by $\bdR$ the boundary at infinity of $\Rm$ (i.e.,\ the set of asymptotic equivalence classes of geodesic rays in $\Rm$).  Another way to view $\bdR$ which will be helpful later is to consider the following equivalence relation on \Hom: $\chi_1 \sim \chi_2$ if and only if $\chi_1 = r\chi_2$ for some $r > 0$.  Define $\bdR$ to be the set of equivalence classes $\{ [\chi] \in$ \Hom $| \hspace{2pt} \chi \neq 0 \}$.  Let $e \in \bdR$ and let $\gamma$ be a geodesic ray defining $e$.  For each $s > 0$, denote by $\hgs$ the half-space in $\Rm$ whose boundary is orthogonal to $\gamma$ with $\hgs \cap \gamma([0,\infty)) = \gamma([s,\infty))$.  Let $\Gamma$ denote the Cayley graph of $G$ with respect to $\mathcal{X}$.  Since the ${\mathbb Z}$-rank of the abelianization of $G$ is $m$, there is an epimorphism $\pi:G \to {\mathbb Z}^m$.  Define $h:\Gamma \to \Rm$ by: $h(g):=\pi(g)$ for all vertices $g \in \Gamma$, and extend linearly on edges.  For each $s \geq 0$, let $\Gamma_{\gamma,s} :=h^{-1}(\hgs)$.  The direction $e \in \Sigma^1(G)$ if and only if for every $s \geq 0$, there exists $\lambda = \lambda(s) \geq 0$ such that any two points $u,v \in \Gamma_{\gamma,s}$ can be joined by a path in $\Gamma_{\gamma,s-\lambda}$ and $s-\lambda(s) \to \infty$ as $s \to \infty$.

In the compactified space $\Rm \cup \bdR$, the compactified half-spaces play
the role of neighborhoods of the point $e \in \bdR$, but this
gives an unsatisfactory topology to $\Rm \cup \bdR$.
From the point of view of topology, it is more natural to have a similar definition to $\Sigma^1(G)$ using ``ordinary" neighborhoods of $e$.  A basis for these neighborhoods consists of ``truncated cones".  For each $s \geq 0$ and each geodesic ray $\gamma$, define the {\it truncated cone} $\cgs := Cone_{\theta}(\gamma) \cap \hgs$ where $Cone_{\theta}(\gamma)$ is the closed cone of angle $\theta$ and vertex $\gamma(0)$, and $\theta := \arctan(\frac{1}{s})$ if $s > 0$ and $\theta := \frac{\pi}{2}$ if $s=0$.  Let $\Delta_{\gamma,s} := h^{-1}(\cgs)$. We say that $e \in \Omega^1(G)$ if and only if there exists $s_0 \geq 0$ such that for each $s \geq s_0$, there exists $\lambda = \lambda(s) \geq 0$ such that any two points $u,v \in \Delta_{\gamma,s}$ can be joined by a path in $\Delta_{\gamma,s-\lambda}$ and $s-\lambda(s) \to \infty$ as $s \to \infty$. We should mention that $\Omega^n(G)$ is always a closed set while $\Sigma^n(G)$ is open. Furthermore, we have $\Sigma^1(G) \supseteq \Sigma^2(G) \supseteq \Sigma^3(G) \supseteq ...$ or $[\Sigma^1(G)]^c \subseteq [\Sigma^2(G)]^c \subseteq [\Sigma^3(G)]^c \subseteq ...$ while we have $\Omega^1(G) \supseteq \Omega^2(G) \supseteq \Omega^3(G) \supseteq ...$.

This paper is organized as follows. In section 2, we review some background on twisted conjugacy classes. Basic properties of the $\Omega^n$-invariants are reviewed in section 3. We prove our main results in section 4 together with applications in section 5.

We thank Daciberg Gon\c calves for providing us the preprint \cite{GK}, Ross Geoghegan for helpful comments on an earlier version of the manuscript, Alexander Fel'shtyn for bringing to our attention the paper \cite{F} in which arguments were sketched for the $R_{\infty}$ property for relatively hyperbolic groups which include free products, and to Enric Ventura for informing us of the work \cite{CMV} where independently the authors claimed that any finite free product of finitely generated freely indecomposable groups has property $R_{\infty}$.

\section{Twisted Conjugacy Classes}

The basic algebraic techniques used in the present paper for showing $R(\varphi)=\infty$ is the relationship among the Reidemeister numbers of group homomorphisms of a short exact sequence. In general, given a commutative diagram of groups and homomorphisms
\begin{equation*}
\begin{CD}
    A    @>{\eta}>>  B  \\
    @V{\psi}VV  @VV{\varphi}V   \\
    A    @>{\eta}>>  B
\end{CD}
\end{equation*}
the homomorphism $\eta$ induces a function $\hat {\eta}:\mathcal R(\psi) \to \mathcal R(\varphi)$ where $\mathcal R(\alpha)$ denotes the set of $\alpha$-twisted conjugacy classes.
Some of the basic facts that will be used throughout this paper are given in the following lemma. For our purposes, we are only concerned with automorphisms while the following result can be proven in greater generality. For more general results, see \cite{daci-peter} and \cite{wong}.

\begin{lem}\label{abelian}
Given an automorphism $\psi:G\to G$ of a finitely generated torsion-free abelian group $G$, $Fix \psi=\{1\}$ iff $R(\psi)<\infty$.
\end{lem}
\begin{proof}
Suppose $G$ is finitely generated torsion-free abelian. Then $G=\mathbb Z^k$ for some positive integer $k$. For any $\varphi: G\to G$, $\# Coker (1-\varphi)<\infty iff |\det (1-\varphi)|\ne 0$ in which case $R(\varphi)=\# Coker (1-\varphi)=|\det (1-\varphi)|$ It follows that $R(\varphi)<\infty$ iff $\varphi$ does not have $1$ as an eigenvalue iff $\varphi(x)=x$ has only trivial solution, i.e., $Fix \varphi=1$.
\end{proof}

\begin{lem}\label{R-facts}
Consider the following commutative diagram
\begin{equation*}\label{general-Reid}
\begin{CD}
    1 @>>> A    @>>>  B @>>>    C @>>> 1 \\
    @.     @V{\varphi'}VV  @V{\varphi}VV   @V{\overline \varphi}VV @.\\
    1 @>>> A    @>>>  B @>>>    C @>>> 1
 \end{CD}
\end{equation*}
where the rows are short exact sequences of groups and the vertical arrows are group automorphisms.
\begin{enumerate}
\item If $R(\overline \varphi)=\infty$ then $R(\varphi)=\infty$.

\item If $|Fix \overline \varphi|<\infty$ and $R(\varphi')=\infty$ then $R(\varphi)=\infty$.

\item If the short exact sequence is a central extension then $R(\varphi)=R(\varphi')R(\bar {\varphi})$.
\end{enumerate}
\end{lem}
\begin{proof}
The homomorphism $p:B\to C$ induces a function $\hat p:\mathcal R(\varphi) \to \mathcal R(\overline{\varphi})$ given by $\hat p([\alpha]_B)=[p(\alpha)]_C$. Since $p$ is surjective, so is $\hat p$. Thus, $(1)$ follows. Similarly, $i:A\to B$ induces a function $\hat i:\mathcal R(\varphi') \to \mathcal R(\varphi)$. Since the sequence is exact, it follows that $\hat i(\mathcal R(\varphi'))=\hat p^{-1}([1]_C)$. The subgroup $Fix \overline{\varphi}$ acts on $\mathcal R(\varphi')$ via $\bar {\theta}\cdot [\alpha]_A=[\theta \alpha \varphi(\theta)^{-1}]_A$ where $\theta \in B$ and $\bar {\theta}\in Fix \overline{\varphi}$. Thus, two classes $[\alpha]_A$ and $[\beta]_A$ are in the same $Fix \overline{\varphi}$-orbit iff $i(\alpha)$ and $i(\beta)$ are in the same Reidemeister class, i.e., $[i(\alpha)]_B=[i(\beta)]_B$. Now, $(2)$ follows immediately. Finally, if the extension is central, $\hat p^{-1}([\bar \alpha]_C)$ is independent of $\bar \alpha$ so that $R(\varphi)=k\cdot R(\varphi')$ and $k$ is the number of distinct Reidemeister classes of $\overline{\varphi}$ in $C$. In other words, $k=R(\overline{\varphi})$ and thus $(3)$ follows.
\end{proof}

\section{The Geometric Invariants}
Let $n$ be a non-negative integer, and let $G$ be a group of type $F_n$ (i.e.,\ $G$ has a $K(G,1)$ complex with a finite $n$-skeleton.).  In this section, we define two invariants of $G$:
	\begin{enumerate}
	\item the Bieri-Neumann-Strebel-Renz (or BNSR) invariants $\SnG$, and
	\item the invariants $\OnG$.
	\end{enumerate}
	
\subsection{The BNSR invariants $\SnG$}\label{sigma}
Let $e \in \bdR$, and let $\gamma$ be a geodesic ray defining $e$.  Associated to $\gamma$ is the function $\beta_{\gamma}:$ \Hom$ \cong \Rm \to {\mathbb R}$ defined by $\beta_{\gamma}(a) := \langle u_e,a-\gamma(0)\rangle / \|u_e\|$ where $\langle\cdot,\cdot\rangle$ is the chosen inner product for $\Rm$, $\| \cdot \|$ is the norm, and $u_e$ is a vector (at 0) pointing toward $e$.  For each $s \in {\mathbb R}$, let $\hgs := \beta^{-1}_{\gamma}([s,\infty))$.  Each $\hgs$ is a closed half-space orthogonal to $\gamma$.

Pick an $n$-dimensional, $(n-1)$-connected CW-complex $X$ on which $G$ acts freely as a group of cell permuting homeomorphisms with $G \backslash X$ a finite complex.  Choose a $G$-map $h:X \to \Rm$ (note that $G$ acts on $\Rm$ by translations), and for each $s \in {\mathbb R}$, denote by $\xgs$ the largest subcomplex of $X$ contained in $h^{-1}(\hgs)$.  Define the {\it BNSR geometric invariants} of $G$, denoted $\SnG$, by $e \in \SnG$ if and only if for every $s \in {\mathbb R}$ and every $-1 \leq p \leq n-1$, there exists $\lambda=\lambda(s) \geq 0$ such that every continuous map $f:S^p \to \xgs$ can be extended to a continuous map $\hat{f}:B^{p+1} \to \xgsl$ and $s-\lambda(s) \to \infty$ as $s \to \infty$.

\subsection{The invariants $\OnG$}\label{omega}
There are invariants analogous to $\SnG$ that replace half-spaces with ``truncated cones".  Let $e \in \bdR$, and let $\gamma$ be a geodesic ray defining $e$.  For each $s \geq 0$, define the {\it truncated cone} $\cgs := Cone_{\theta}(\gamma) \cap \hgs$ where:
	\begin{enumerate}
	\item $\theta := \arctan(\frac{1}{s})$ if $s > 0$ and $\theta := \frac{\pi}{2}$ if $s=0$, and
	\item $Cone_{\theta}(\gamma)$ is the closed cone of angle $\theta$ and vertex $\gamma(0)$.
	\end{enumerate}
Choose $X$ and $h$ as before.
Denote by $\ygs$ the largest subcomplex of $X$ contained in $h^{-1}(\cgs)$.
Define $\OnG$ by $e \in \OnG$ if and only if there exists $s_0 \geq 0$ such that for every $s \geq s_0$ and each $-1 \leq p \leq n-1$, there exists $\lambda=\lambda(s) \geq 0$ such that every continuous map $f:S^p \to \ygs$ can be extended to a continuous map $\hat{f}:B^{p+1} \to \ygsl$ and $s-\lambda(s) \to \infty$ as $s \to \infty$.

The following theorem relates the invariants $\SnG$ and $\OnG$.
\begin{thm}\label{hemi}\cite[Theorem 3.1]{K1}
Let $e \in \bdR$.  Then $e \in \OnG$ if and only if $\ep \in \SnG$ for every $\ep$ in an open $\frac{\pi}{2}$-neighborhood of $e$.
\end{thm}

Given $\Sigma^n(G)$, we can completely determine
$\Omega^n(G)$: for each $e \in \bdR$, $e \in \Omega^n(G)$ if and only if
the open $\frac{\pi}{2}$-neighborhood of $e$ is in $\Sigma^n(G)$.
However, it is not the case that $\Omega^n(G)$ completely determines
$\Sigma^n(G)$; examples of such groups are given in \cite[\S~1.3]{K1}.

The following theorem completely describes $\OnGH$ in terms of $\OnG$ and $\OnH$.  This theorem will be useful in \S~\ref{app}.
\begin{thm}\label{product-formula}\cite[Theorem 3.8]{K2}
$\OnGH = \OnG \circledast \OnH$ where $\circledast$ represents the spherical join.
\end{thm}

It should be noted that no such result exists for $\Sigma^n(G \times H)$.  In \cite{schutz}, a counterexample to the product conjecture for $\Sigma^n$ was given.  However, in \cite{BG}, the product conjecture for the homological version of the Bieri-Neumann-Strebel-Renz invariants is proven over a field.

\section{The Main Result}

In this section, we obtain information about $R_{\infty}$ when $\Omega^n(G)$ is a non-empty finite set.

\begin{prop}\label{finite12}
Let $n$ be a positive integer, and let $G$ be a group of type $F_n$ . If $0<\# \Omega^n(G) < \infty$ then $\# \Omega^n(G) \in \{1,2\}$.  Moreover, if $\# \Omega^n(G) = 2$ then the two points in $\Omega^n(G)$ are antipodal.
\end{prop}
\begin{proof} Suppose $e,\ep \in \Omega^n(G)$, and suppose the distance between $e$ and $\ep$ is less that $\pi$.  Let $\epp$ be a point on the geodesic joining $e$ and $\ep$.  Then the $\frac{\pi}{2}$-neighborhood of $\epp$ would be contained in the union of the $\frac{\pi}{2}$-neighborhoods of $e$ and $\ep$.  By Theorem~\ref{hemi}, since $e, \ep \in \Omega^n(G)$, the $\frac{\pi}{2}$-neighborhoods of $e$ and $\ep$ are contained in $\Sigma^n(G)$.  Thus, the $\frac{\pi}{2}$-neighborhood of $\epp$ is contained in $\Sigma^n(G)$, so by Theorem~\ref{hemi}, $\epp \in \Omega^n(G)$.  Therefore, if $\Omega^n(G)$ is finite and $e \in \Omega^n(G)$, no other points within a distance of $\pi$ of $e$ can be in $\Omega^n(G)$.  The only other possible point would be the point antipodal to $e$.
\end{proof}

Next, we show that $\OnG$ is invariant under automorphisms of $G$.  Although this result is true for $\SnG$ and we use similar techniques in the proof for $\OnG$, we feel it necessary to prove the following proposition since $\SnG$ cannot be completely determined from $\OnG$.  Thus, the result does not necessarily follow from the fact that $\SnG$ is invariant under automorphisms.

\begin{prop}\label{invariant}
Let $n$ be a positive integer, and let $G$ be a group of type $F_n$.  Then $\OnG$ is invariant under automorphisms of $G$.
\end{prop}

\begin{proof} Let $[\chi] \in \OnG$, and let $\varphi \in {\rm Aut}(G)$.  Define $\tilde{\varphi}: \Hom \to \Hom$ by $\tilde{\varphi}(\alpha) = \alpha \circ \varphi$ for any $\alpha \in \Hom$.  We will show by induction that $[\tilde{\varphi}(\chi)] \in \OnG$.  Suppose $G$ has generating set $\mathcal{S}$, and let $\Gamma(\mathcal{S})$ be the Cayley graph of $G$ with respect to $\mathcal{S}$.  Then the set $\mathcal{S}^{\prime} := \varphi(\mathcal{S})$ generates $G$, and let $\Gamma(\mathcal{S}^{\prime})$ be the Cayley graph of $G$ with respect to $\mathcal{S}^{\prime}$.  The induced map $\varphi^{\prime}: \Gamma(\mathcal{S}) \to \Gamma(\mathcal{S}^{\prime})$ is a cellular map.  If $h: \Gamma(\mathcal{S}) \to \Hom$ is a $G$-map, then $\bar{h}: \Gamma(\mathcal{S})^{\prime} \to \Hom$ defined by $\bar{h} := \tilde{\varphi} \circ h \circ \varphi^{\prime -1}$ is a $G$-map.  Let $\gamma$ be a geodesic ray in the direction of $[\tilde{\varphi}(\chi)]$, so $\gamma^{\prime} = \tilde{\varphi}^{-1}(\gamma)$ is a geodesic ray in the direction of $[\chi]$.  Let $C_{\gamma,s}$ be a truncated cone, so $C_{\gamma^{\prime},s} = \tilde{\varphi}^{-1}(C_{\gamma,s})$ is a truncated cone, and $\bar{h}^{-1}(C_{\gamma,s}) = \varphi^{\prime}(h^{-1}(C_{\gamma^{\prime},s}))$.  Since $[\chi] \in \Omega^1(G)$, we know that $h^{-1}(C_{\gamma^{\prime},s})$ satisfies the required connectivity property, and since $\varphi^{\prime}$ is a cellular map (and a path in the Cayley graph is a series of generators), we have that $\bar{h}^{-1}(C_{\gamma,s})$ satisfies the required connectivity property.  By induction on $1 \leq k \leq n$, the map $\varphi^{\prime}: X^k \to X^k$ is a cellular map, and we get that $[\tilde{\varphi}(\chi)] \in \OnG$.
\end{proof}

A point $[\chi]\in \partial_{\infty}\mathbb R^m$ (boundary of $\Hom$) is {\it rational} (or discrete as in \cite{GK}) if $\chi(G)$ is infinite cyclic. Here is our main theorem.

\begin{thm}\label{main}
Let $n$ be a positive integer, and let $G$ be a group of type $F_n$ with $0<\# \Omega^n(G) < \infty$.  Suppose that $\OnG$ contains only rational points.
\begin{enumerate}
\item If $\# \Omega^n(G)=1$, then $G$ has property $R_{\infty}$.
\item If $\# \Omega^n(G)=2$, then there exists a normal subgroup $N \lhd {\rm Aut}(G)$ with $[{\rm Aut}(G):N]=2$ such that $R(\varphi)=\infty$ for every $\varphi \in N$.
\end{enumerate}
\end{thm}

\begin{proof} Suppose $\# \Omega^n(G) = 1$ and $[\chi] \in \OnG$.  Let $N = ker(\chi)$, let $V :=$ Hom$_{\mathbb{Z}}(G/N,\mathbb{R})$, and let $\varphi \in {\rm Aut}(G)$.  {Since $[\chi]$ is rational, we have $G/N$ has rank 1, so $V$ is 1-dimensional.} From Proposition \ref{invariant}, $[\tilde{\varphi}(\chi)] \in \OnG$ so $[\tilde{\varphi}(\chi)] =[\chi]$ implies that $\chi \circ \varphi =c \chi$ for some $c\in \mathbb Z$ and hence $\varphi (N)\subseteq N$, so $N$ is characteristic in $G$. The automorphism $\varphi$ induces the following two maps: $\overline{\varphi}: G/N \to G/N$ defined by $\overline{\varphi}(gN) = \varphi(g)N$, and $\hat{\varphi}: V \to V$ defined by $\hat{\varphi}(\alpha)(gN) = \alpha(\varphi(g)N)$.  We will show that $\overline{\varphi}$ is the identity map.  Since $\varphi$ is invertible and $N$ is characteristic, we have that $\hat{\varphi}$ is invertible and that $\{ [\overline{\chi}] \}$ is a basis for $V$ where $\overline{\chi}:G/N \to \mathbb R$ is the induced map from $\chi$.  Since $\hat{\varphi}$ is invertible, $c = \pm 1$, but $c \neq -1$ since $[-\chi] \not\in \OnG$.  Therefore, $\hat{\varphi}(\overline{\chi}) = \overline{\chi}$.  Thus, we have $\hat{\varphi}(\overline{\chi})(gN) = \overline{\chi}(\overline{\varphi}(gN))$ which implies $\overline{\chi}(gN) = \overline{\chi}(\varphi(g)N)$.  Therefore, $g^{-1}\varphi(g) \in N$, so $\varphi(g)N = gN$ which implies $\overline{\varphi}(gN) = \varphi(g)N = gN$. Now $G/N$ is a finitely generated abelian group of rank $1$ so that $(G/N)/\{\text{torsion}\} \cong \mathbb Z$. Moreover, $\overline \varphi$ is the identity on $G/N$ so it also induces the identity on $(G/N)/\{\text{torsion}\}$. It follows from Lemma \ref{abelian} that $R(\overline{\varphi})=\infty$ so that $G$ has the property $R_{\infty}$ by Lemma \ref{R-facts}.

Now suppose $\# \OnG = 2$.  By Proposition~\ref{finite12}, the two points must be antipodal, say $[\chi]$ and $[-\chi]$.  Let $H$ be the subgroup of automorphisms $\varphi$ of $G$ that fix $[\chi]$ under the induced map $\tilde \varphi$ as in the proof of Proposition \ref{invariant}.  Then $H$ is an index two normal subgroup of ${\rm Aut}(G)$, and by the preceding paragraph, every automorphism in $H$ has an infinite number of twisted conjugacy classes.
\end{proof}

Now we give a simpler and shorter proof of the fact that the solvable Baumslag-Solitar groups $BS(1,n)$ for $n\ge 2$ have property $R_{\infty}$ \cite{fel-daci}.

\begin{cor}
$BS(1,n)$ has property $R_{\infty}$.
\end{cor}
\begin{proof} A similar argument as in \cite[Example]{K2} for $BS(1,2)$ shows that $\Omega^1(BS(1,n))=\{+\infty\}$ consisting of one rational point. The assertion follows immediately from Theorem \ref{main}.
\end{proof}

Next, we compare our results with those of \cite{GK} which uses the $\Sigma^n$-theory.
The following results are due to D. Gon\c calves and D. Kochloukova \cite{GK}.

\begin{thm}\label{gk1}\cite[Corollary 6]{GK}
Let $G$ be a finitely generated group such that $\Sigma^1(G)^c=\{[\chi_1],...,[\chi_m]\}$ is a non-empty finite set of discrete points. Then ${\rm Aut}(G)$ has a finite index subgroup $H$ such that $R(\varphi)=\infty$ for all $\varphi \in H$.
\end{thm}

\begin{thm}\label{gk2}\cite[Theorem 4]{GK}
Let $G$ be as in Theorem \ref{gk1} and let $N=\bigcap_i {\rm Ker} \chi_i$. If the images of $\{\chi_1,...,\chi_m\}$ in $V={\rm Hom}_{\mathbb Z}(G/N,\mathbb R)$ form a basis for $V$ as a real vector space, then $G$ has property $R_{\infty}$.
\end{thm}

Recall \cite{BrGu} the generalized Thompson's groups $F_{n,0}$ for $n \geq 2$ are defined with the following presentation
\begin{equation*}
\langle x_0,x_1,x_2,\ldots | \hspace{3pt} x_j^{-1}x_ix_j = x_{i+n-1} \hspace{3pt} {\rm for} \hspace{3pt} 0 \leq j < i \rangle.
\end{equation*}

When $n=2$, $F_{2,0}$ is the original Thompson group $F$ which is known to have property $R_{\infty}$ \cite{BFG}.
In \cite[Corollary 9]{GK}, it was shown that the groups $F_{n,0}$ have property $R_{\infty}$ by using the fact $\Sigma^1(F_{n,0})^c=\{[\rho],[\lambda]\}$
so that Theorem \ref{gk2} is applicable. Furthermore, using Theorem \ref{gk1}, Gon\c calves and Kochloukova showed that if $G$ is either a finitely generated metabelian group of finite Pr\"ufer rank \cite[Lemma 11]{GK} or the pure symmetric automorphism group of free groups of finite rank \cite[Theorem 18]{GK} or a one-ended non-abelian limit group that is not a surface group \cite[Theorem 22]{GK}, then ${\rm Aut}(G)$ contains a finite index subgroup $H$ such that $R(\varphi)=\infty$ for every $\varphi \in H$. It follows from \cite{BGK} that all $\Sigma^n(F_{n,0})$ are infinite and independent of $n$ and hence $\Omega^n(F_{n,0})$ is infinite for all $n$.
Thus, the $R_{\infty}$ for $F_{n,0}$ cannot be deduced from our result Theorem \ref{main}. On the other hand, the following examples serve as applications of the $\Omega^n$-theory while the $\Sigma^n$-theory method of Gon\c calves and Kochloukova \cite{GK} is not applicable.

\begin{ex}\label{ex4-1}
Let $G=BS(1,2)\times \mathcal F_n$ where $\mathcal F_n$ is the free group of rank $n\ge 2$. Then, using the product formula of Theorem \ref{product-formula} for $\Omega^1$, we have $\Omega^1(G)=\{+\infty\}$ consisting of one rational point and hence $G$ has property $R_{\infty}$. On the other hand,
\begin{equation*}
    \begin{aligned}
          \left[\Sigma^1(BS(1,2)\times \mathcal F_n)\right]^c&=[\Sigma^1(BS(1,2))]^c \circledast [\Sigma^0(\mathcal F_n)]^c \cup [\Sigma^0(BS(1,2))]^c \circledast [\Sigma^1(\mathcal F_n)]^c \\
                                        &=(\{-\infty\}\circledast \emptyset) \cup (\emptyset \circledast \emptyset^c) =\{-\infty\} \cup S^{n-1}
    \end{aligned}
\end{equation*}
so that $[\Sigma^1(BS(1,2)\times \mathcal F_n)]^c$ is not finite and hence Theorem \ref{gk2} is not applicable.
\end{ex}

\begin{ex}\label{ex4-2}
Let $G=\mathcal F_n \times \mathbb Z$, for $n\ge 2$. Then $\Omega^1(G)=\{\pm \infty\}$ consisting of two antipodal rational points while
\begin{equation*}
    \begin{aligned}
          \left[\Sigma^1(G)\right]^c&=([\Sigma^1(\mathcal F_n)]^c \circledast [\Sigma^0(\mathbb Z)]^c ) \cup ([\Sigma^0(\mathcal F_n)]^c \circledast [\Sigma^1(\mathbb Z)]^c) \\
               &=(S^1 \circledast \emptyset) \cup (\emptyset \circledast \emptyset) = S^{n-1}.
    \end{aligned}
\end{equation*}
Thus, it follows from Theorem \ref{main} that ${\rm Aut}(\mathcal F_n \times \mathbb Z)$ has an index $2$ subgroup $H$ such that $R(\varphi)=\infty$ for all $\varphi \in H$. On the other hand, Theorem \ref{gk1} does not apply.
In fact, $\mathcal F_n \times \mathbb Z$ has property $R_{\infty}$ since it is a GBS (Generalized Baumslag-Solitar) group  \cite{levitt}. Furthermore, it is known \cite{TW2} that
any group quasi-isometric to $\mathcal F_n \times \mathbb Z$ has $R_{\infty}$.
\end{ex}

\begin{ex}\label{ex4-3}
Let $B_n$ be the Artin braid group on $n$ strands and $G=\pi_1(K)$, the Klein bottle group.  Since the commutator subgroup of $B_n$ is finitely generated \cite{GL}, we have $\Omega^1(B_n) = \{\pm \infty\}$ while $[\Sigma^1(B_n)]^c=\emptyset$, and the same is true for the Klein bottle group, $\Omega^1(G)=\{\pm \infty\}$ while $[\Sigma^1(G)]^c=\emptyset$.
\end{ex}

\section{Applications}\label{app}

In this section, we make use of the product formula for $\Omega^n$ as in Theorem \ref{product-formula} to produce new examples of groups with property $R_{\infty}$ and of groups of which the automorphism group contains an index 2 subgroup with property $R_{\infty}$.

For any positive integer $k$, $i\in \{0,1,2\}$, we denote by $\mathcal O^k_i$ the family of finitely generated groups $G$ such that $\# \Omega^k(G)=i$ with all points in $\Omega^k(G)$ rational. Since $\Omega^1(G) \supseteq \Omega^2(G) \supseteq ...$, $\mathcal O^n_i$ is a sub-family of $\mathcal O^m_j$ for any $m\le n$ where $j\le i$.

\begin{thm}
Let $m \le n$ be positive integers and $K\in \mathcal O^m_0$.
\begin{enumerate}
\item If $H\in \mathcal O^n_1$ then $H\times K$ has property $R_{\infty}$.
\item If $H\in \mathcal O^n_2$ then ${\rm Aut}(H\times K)$ contains a subgroup $N$ of index $2$ such that $R(\varphi)=\infty$ for every $\varphi \in N$.
\end{enumerate}
\end{thm}

\begin{proof}
(1) Since $K\in \mathcal O^m_0$ and $m\le n$, we have $\emptyset=\Omega^m(K) \supseteq \Omega^n(K)$. It follows from Theorem \ref{product-formula} that $\#\Omega^n(H\times K)=\#(\Omega^n(H)\circledast \Omega^n(K))=\#\Omega^n(H)=1$ since $H\in \mathcal O^n_1$. By (1) of Theorem \ref{main}, $H\times K$ has property $R_{\infty}$.

(2) Now, $H\in \mathcal O^n_2$ implies that $\#\Omega^n(H\times K)=\#\Omega^n(H)=2$ and thus the assertion follows from (2) of Thoerem \ref{main}.
\end{proof}

Next, we extend our results to free products of groups.

\begin{thm}\label{free-product}
Let $G=\ast_{i=1}^n A_i$ be a finite free product of $n$ non-trivial freely indecomposable finitely generated groups, $n\ge 2$. If one of the following conditions is satisfied then $G$ has the property $R_{\infty}$.
\begin{enumerate}
\item Each $A_i$ is finite.
\item There exists $j, 1\le j\le n$ such that $A_j\in \mathcal O^m_1$ and for $i\ne j$, $A_i\in \mathcal O^{k_i}_0$ with $k_i\le m$.
\item The direct product $\bar G=\prod A_i$ has property $R_{\infty}$ and $A_i$ is abelian and non-isomorphic to $\mathbb Z$ for some $i$.
\end{enumerate}
\end{thm}
\begin{proof} Let $C(G)$ be the Cartesian subgroup of $G$, that is, the kernel of the canonical map $\ast_{i=1}^n A_i \to \prod A_i=\bar G$.

(1) If each $A_i$ is finite, none of them is isomorphic to $\mathbb Z$. It follows from \cite{collins} that $C(G)$ is free and characteristic. Moreover, since each $A_i$ is finite, $C(G)$ is finitely generated and thus is a free group of finite rank $r$. If $G=\mathbb Z_2 \ast \mathbb Z_2=D_{\infty}$, it has already been shown in \cite{daci-peter3} that the infinite dihedral group has the property $R_{\infty}$. Otherwise, the free group $C(G)$ has rank $r\ge 2$ and hence is a finitely generated non-elementary Gromov hyperbolic group. It follows from \cite{levitt-lustig} that $C(G)$ has property $R_{\infty}$. Since $\bar G$ is finite, it follows from (2) of Lemma \ref{R-facts} that $G$ has property $R_{\infty}$.

(2) Since $A_j\in \mathcal O^m_1$, $A_j$ is not isomorphic to $\mathbb Z$ for $\#\Omega^n(\mathbb Z)=2$. Similarly, $A_i\ncong \mathbb Z$ for $A_i\in \mathcal O^{k_i}_0$. Therefore, it follows from \cite{collins} that $C(G)$ is characteristic. Since $\bar G$ has $R_{\infty}$, it follows from (1) of Lemma \ref{R-facts} that $G$ has property $R_{\infty}$. Note that $C(G)$ need not be finitely generated so the result of \cite{levitt-lustig} is not applicable.

(3) It follows from \cite{NT} that $C(G)$ is characteristic. Then the result follows from (1) of Lemma \ref{R-facts}.
\end{proof}

\begin{ex}\label{ex5-2}
Take $A_1=BS(1,m)$ and for $2\le i\le n$, $A_i$ to be a non-trivial finite group. Then, $A_1 \in O^1_1$ and $A_i \in O^1_0$ for $i\ne 1$. By (2) of Theorem \ref{free-product}, $G=BS(1,m) \ast A_2 \ast ... \ast A_n$ has property $R_{\infty}$.
\end{ex}

\begin{ex}\label{ex5-3}
Let $A_1=\pi_1(K)$ be the fundamental group of the Klein bottle and $A_i=\mathbb Z$ for $i=2,..., n-1$. By Theorem 2.4 of \cite{daci-peter3}, the group $G'=\prod_{i=1}^{n-1} A_i$ has property $R_{\infty}$. Now, let $A_n$ be a non-trivial finite abelian group. Since $G'$ is torsion-free and $A_n$ is finite, it follows that $G'$ is characteristic in the direct product $G=\prod_{i=1}^n A_i=G' \times A_n$. Thus (2) of Lemma \ref{general-Reid} implies that $G$ has property $R_{\infty}$. Now, (3) of Theorem \ref{free-product} shows that the group $\ast_{i=1}^n A_i=\pi_1(K) \ast \mathcal F_{n-1} \ast A_n$ has property $R_{\infty}$.
\end{ex}

Theorem \ref{free-product} shows that many free products of freely indecomposable groups possess property $R_{\infty}$ while the $\Omega^n$ invariants of such groups are necessarily empty. Furthermore, in contrast to Theorem \ref{main} where groups with $\#\Omega^n(G)=1$ must have property $R_{\infty}$, there are examples of groups $G$ with property $R_{\infty}$ that have $\#\Omega^n(G)=2$ or $\#\Omega^n(G)=\infty$.

\begin{ex}\label{ex5-4}
Here, we list some known examples of groups $G$ with property $R_{\infty}$ and various values for the cardinality of $\Omega^n(G)$.
\begin{enumerate}
\item $(\#\Omega^n(G)=0$ and $R_{\infty})$: Groups satisfying one of the conditions in Theorem \ref{free-product}; free groups $\mathcal F_n$ with $n\ge 2$; lamplighter groups $L_n$ with $gcd(n,6)>1$.
\item $(\#\Omega^n(G)=1$ and $R_{\infty})$: $BS(1,n)$.
\item $(\#\Omega^n(G)=2$ and $R_{\infty})$: $\mathcal F_n \times \mathbb Z$ with $n\ge 2$; Artin Braid group on $3$ strands $B_3$; $\pi_1(K)$; certain right angled Artin groups.
\item $(\#\Omega^n(G)=\infty$ and $R_{\infty})$: Thompson's group $F_{n,0}$; $\pi_1(K) \times \mathbb Z^n$ with $n\ge 1$.
\end{enumerate}
\end{ex}

\noindent
{\Small
\begin{tabbing}
Nic Koban \hphantom{xxxxxxxxxxxxxxxxxxxxxxxxxx} \= Peter Wong\\
Department of Mathematics \> Department of Mathematics\\
University of Maine Farmington \> Bates College\\
Farmington, ME 04938 \> Lewiston, ME 04240 \\
USA \> USA\\
nicholas.koban@maine.edu \> pwong@bates.edu\\
\end{tabbing}
}
\end{document}